\documentclass[12pt]{article}

\usepackage{amssymb}
\usepackage{amsmath}

\usepackage{amsbsy}

\begin{document}%



\title{\bf Exact expansions of Hankel transforms and related integrals}

\author{
A.V. Kisselev\thanks{Electronic address:
alexandre.kisselev@ihep.ru} \\
{\small A.A.~Logunov Institute for High Energy Physics, NRC
``Kurchatov Institute''}, \\
{\small 142281 Protvino, Russia} }

\date{}

\maketitle

\begin{abstract}
The Hankel transform $\mathcal{H}_n [f(x)](q)=\int_0^{\infty} \!\!
\, x f(x) J_n(q x) dx$ is studied for integer $n\geqslant -1$ and
positive parameter $q$. It is proved that the Hankel transform is
given by uniformly and absolutely convergent series in reciprocal
powers of $q$, provided special conditions on the function $f(x)$
and its derivatives are imposed. It is necessary to underline that
similar formulas obtained previously are in fact asymptotic
expansions only valid when $q$ tends to infinity. If one of the
conditions is violated, our series become asymptotic series. The
validity of the formulas is illustrated by a number of examples.
\end{abstract}



\section{Introduction} %
\label{sec:int}

The calculation of the Hankel transform
\begin{equation}\label{Hnu_def}
\mathcal{H}_\nu [f(x)](q) = \int\limits_0^{\infty} \!\!x f(x)
J_\nu(q x) \,dx
\end{equation}
is a problem which arises in mathematical physics (see, for
instance, \cite{Tranter}, \cite{Debnath}), as well as in high energy
nuclear and particle physics \cite{Moliere}-\cite{Collins}. In
particular, an asymptotic expansion of $\mathcal{H}_\nu [f(x)](q)$
is often needed, as $q \rightarrow \infty$. For the first time
asymptotic representations for $\mathcal{H}_0[f(x)](q)$ and
$\mathcal{H}_1[f(x)](q)$ were given by Willis \cite{Willis} (see
also \cite{Tranter}) without investigating convergence conditions.

During the years that followed, asymptotics of $\mathcal{H}_\nu
[f(x)](q)$ and of more general transforms were derived by a number
of authors using different conditions imposed on the function $f(x)$
\cite{Slonovskii}-\cite{Galapon}.%
\footnote{Methods of evaluating asymptotics of other integrals can
be found in \cite{Bleistein}.}
The aim of all these papers was to obtain \emph{an asymptotic
expansion} of $\mathcal{H}_\nu [f(x)](q)$ valid when $q$ tends to
infinity.

The goal of the present paper is to find such conditions on $f(x)$
and its derivatives under which $\mathcal{H}_\nu [f(x)](q)$ with
integer $\nu$ is presented by \emph{an uniformly and absolutely
convergent series} in reciprocal powers of  $q$, for all $q>a>0$. We
start from the Hankel transforms of order zero, see
section~\ref{sec:zero}. The Hankel transforms of order one and two
are studied in section~\ref{sec:one}. In section~\ref{sec:nu} a
generalization for the Hankel transform of order $n \geqslant -1$ is
presented.

\section{Ramanujan's master theorem} %
\label{sec:Ramanujan}

Our further consideration will be based on the Ramanujan master
theorem (see its proof provided by Hardy in \cite{Hardy}). We use
the version of this theorem presented in \cite{Amdeberham}.

\textsc{Ramanujan's master theorem}. Let $\varphi(z)$ be analytic
(single-valued) function, defined on a half-plane
\begin{equation}\label{Ram_h_half_plane}
H(\delta) = \{z \in \mathrm{\mathbb{C}}\!:  \mathfrak{R}(z)
\geqslant -\delta \}
\end{equation}
for some $0 < \delta <1$. Suppose that, for some $A < \pi$,
$\varphi(z)$ satisfies the growth condition
\begin{equation}\label{Ram_phi_cond}
|\varphi(\sigma + i\tau)| < C e^{P\sigma + A|\tau|}
\end{equation}
for all $z = \sigma + i\tau \in H(\delta)$. Then the identity
\begin{equation}\label{master_theorem}
\int\limits_0^\infty x^{s-1} \left[ \varphi(0) - x \varphi(1) + x^2
\varphi(2) - \ldots \right] dx = \frac{\pi}{\sin (\pi s)} \, \varphi
(-s)
\end{equation}
holds for all $0 < \mathfrak{R}(s) < \delta$.%
\footnote{The extension of Ramanujan's master theorem for $A \in
(0,\pi]$ was done in \cite{Claudhry}.}

Let us define the function $\psi(s) = \varphi (-s)/\sin (\pi s)$.
Due to the growth condition imposed on $\varphi (s)$, $\lim_{|t|
\rightarrow \infty} \psi(a+it) = 0$ uniformly in $a \in
[\alpha,\beta]$ for all $[\alpha,\beta] \in (0,\delta)$. Using
inverse Mellin transform, we get from the Ramanujan master theorem
($0 < x < e^{-P}$) \cite{Claudhry}
\begin{equation}\label{Ramanujan_Mellin}
\sum_{m=0}^{\infty} (-x)^m \varphi(m) = \frac{1}{2 i}
\int\limits_{a-i\infty}^{a+i\infty} \frac{\varphi(-s)}{\sin (\pi s)}
\,x^{-s} \,ds = -\frac{1}{2 i} \int\limits_{-a-i\infty}^{-a+i\infty}
\frac{\varphi(s)}{\sin (\pi s)} \,x^{s} \,ds \;,
\end{equation}
where $0 < a < \delta$. Observe that $1/\sin(\pi s)$ has simple
poles at $s=-n$ for $n=0,1,2, \ldots$ with the residue $(-1)^n$.
Then the Cauchy residue theorem yields \eqref{Ramanujan_Mellin}.

\textsc{Corollary} Let $\varphi(z)$ be meromorphic function on the
half-plane \eqref{Ram_h_half_plane} which obeys the conditions
\eqref{Ram_phi_cond} of the Ramanujan master theorem. Suppose that
$\varphi(z)$  has a finite number of poles at the points $z_i$
($i=1,2, \ldots N$) with $0> \mathfrak{Re}(z_i) > -\delta$. Then
\begin{equation}\label{Ramanujan_Mellin_gen}
\sum_{m=0}^{\infty} (-x)^m \varphi(m) + \pi\sum_{i=1}^N \mathrm{Res}
\!\left[\frac{\varphi(-s) \,x^{-s}}{\sin (\pi s)}; -z_i  \right] =
\frac{1}{2 i} \int\limits_{a-i\infty}^{a+i\infty}
\frac{\varphi(-s)}{\sin (\pi s)} \,x^{-s} \,ds \;,
\end{equation}
where $-\min[\mathfrak{Re}(z_1), \ldots, \mathfrak{Re}(z_N)] < a <
\delta$.

\section{Hankel transform of order zero} %
\label{sec:zero}

Consider the Hankel transform of order zero
\begin{equation}\label{H0_def}
\mathcal{H}_0 [f(x)](q) = \int\limits_0^{\infty} \!\!x f(x) J_0(q x)
\,dx \;,
\end{equation}
where $q$ is a positive parameter. We don't assume that $q$ is
large.

\textsc{Theorem 1}. The Hankel transform \eqref{H0_def} can be
presented by absolutely and uniformly convergent series for $q>a>0$
\begin{align}\label{H0_final}
\mathcal{H}_0 [f(x)](q) &= \frac{1}{q^3} \sum_{m=0}^\infty
(-1)^{m+1} \,\frac{\Gamma (2m+2)}{\Gamma^2 (m+1)} \,f^{(2m+1)}(0) \,
(2q)^{-2m} \nonumber \\
&= \frac{4}{\pi q^3} \sum_{m=0}^\infty (-1)^{m+1}
\,\frac{[\Gamma(m+3/2)]^2}{\Gamma(2m+2)} \,f^{(2m+1)}(0) \left(
\frac{q}{2} \right)^{-2m} \nonumber \\
&= \frac{2}{q^2} \sum_{k=0}^\infty \,\frac{\Gamma (k/2 + 1)}{\Gamma
(- k/2) \,\Gamma(k+1)} \,f^{(k)}(0) \left( \frac{q}{2}
\right)^{\!\!-k} ,
\end{align}
provided that
\begin{description}
  \item[(1)]
$f(-x) \neq f(x)$;
  \item[(2)]
$f(x)$ is a regular function at a point $x=0$, and its Taylor series
has the form
\begin{equation}\label{g_Taylor_series}
f(x) = \sum_{k=0}^\infty \frac{g^{(k)}(0)}{k!} \,(-x)^k \;,
\end{equation}
where $(-1)^k g^{(k)}(x)$ is the $k$-th derivative of $f(x)$;
  \item[(3)]
$x^{1/2} f(x) \rightarrow 0$, as $x \rightarrow \infty$;
  \item[(4)]
$g^{(s)}(0)$ is a regular (analytic single-valued) function defined
on a half-plane
\begin{equation}\label{half_plane}
H(\delta) = \{s\in \mathrm{\mathbb{C}}\!: \mathfrak{R}(s) \geqslant
-\delta \}
\end{equation}
for some $3/4 < \delta < 1$ and satisfies the growth condition
\begin{equation}\label{f_growth}
|f^{(s)}(0)| < C a^{\mathfrak{R}(s)}
\end{equation}
for some $a>0$,  and all $s \in H(\delta)$.
\end{description}

\textsc{Proof.} Let us prove that our series \eqref{H0_final} is
equal to the Hankel transform \eqref{H0_def}. The r.h.s. of
eq.~\eqref{H0_final} can be rewritten as
\begin{equation}\label{S0_def}
S_0(q) = \frac{1}{q^3} \sum_{m=0}^\infty (-1)^{m} \,\frac{\Gamma
(2m+2)}{\Gamma^2 (m+1)} \,g^{(2m+1)}(0) \, (2q)^{-2m} \;.
\end{equation}
Since
\begin{equation}\label{Gamma_ratio_asymp_1}
\frac{\Gamma (2m+2)}{\Gamma^2 (m+1)}\bigg|_{m\rightarrow \infty} =
2^{\,2m+1} \sqrt{\frac{m}{\pi}}  \left[ 1 + \mathrm{O}(m^{-1})
\right] ,
\end{equation}
the function
\begin{equation}\label{phi_H0}
\varphi(s) = \frac{\Gamma (2s+2)}{\Gamma^2 (s+1)} \, g^{(2s+1)}(0)
\end{equation}
obeys all the condition of the Ramanujan master theorem for $q > a$.
Using formula \eqref{Ramanujan_Mellin}, we obtain
\begin{equation}\label{S0_contour_initial}
S_0(q) = 4i \int\limits_{c-i\infty}^{c+i\infty} \frac{1}{\sin(\pi
s)} \,\frac{\Gamma(2s+2)}{\Gamma^2(s+1)} \,g^{(2s+1)}(0)
\,(2q)^{-2s-3} ds \;,
\end{equation}
where $-\delta < c < -3/4$. Taking into account that $\Gamma(2s+2) =
2^{2s+1}\Gamma(s+1) \Gamma(s+3/2)/\sqrt{\pi}$, and $\sin(\pi
s)\Gamma(s+1) = -\pi/\Gamma(-s)$ \cite{Bateman_vol_1}, we find
\begin{align}\label{S0_contour_median}
S_0(q) &= \frac{1}{i \pi^{3/2}} \int\limits_{c-i\infty}^{c+i\infty}
\Gamma(s+3/2) \Gamma(-s) \,g^{(2s+1)}(0) \,q^{-2s-3} \,ds \nonumber
\\
&= \frac{1}{4i\pi} \int\limits_{c'-i\infty}^{c'+i\infty}
\frac{\Gamma(s'/2 + 1) \Gamma(-s')}{\Gamma(-s'/2)} \,g^{(s')}(0)
\,\left( \frac{q^2}{4} \right)^{\!\!-(s'/2+1)} ds' \;,
\end{align}
were $-1 < 1 - 2\delta  < c' < -1/2$, and the change of variables
$s'=2s+1$ is made. According to \cite[2.12.2.2]{Bateman_vol_2},
\begin{equation}\label{Bessel_integral}
\frac{\Gamma(s'/2+1)}{\Gamma(-s'/2)} \, \left( \frac{q^2}{4}
\right)^{\!\!-(s'/2+1)} = 2 \!\int\limits_0^\infty \!\!x^{s'+1}
J_0(x q) \,dx  \;.
\end{equation}
Integral in \eqref{Bessel_integral} converges, since $-1 <
\mathfrak{R}(s') < -1/2$. As a result, we get
\begin{equation}\label{S0_contour_final}
S_0(q) = \frac{1}{2\pi i} \int\limits_{c'-i\infty}^{c'+i\infty}
\Gamma(-s') \,g^{(s')}(0)  \int\limits_0^\infty x^{s'+1} J_0(x q)
\,dx \;.
\end{equation}
The Cauchy residue theorem with the relation $(-1)^n g^{(n)}(0) =
f^{(n)}(0)$ yields eq.~\eqref{H0_final}. Q.E.D.

Let us emphasize the importance of the condition \textbf{(1)} of
theorem~1 which says that the even-order and odd-order derivatives
of $f(x)$
have opposite sings.%
\footnote{All even-order derivative are positive, while all
odd-order derivatives are negative, or vice versa.}

The r.h.s. in \eqref{H0_final} is \emph{the uniformly and absolutely
convergent series} for all $q > a$. Consequently, it is the exact
representation of the Hankel transform \eqref{H0_def}, contrary to
the results of papers \cite{Slonovskii}-\cite{Soni} in which similar
expansions were derived as \emph{asymptotic series} for large $q$.

Let us see how does our formula \eqref{H0_final} work? It is
worthwhile considering several examples.\\
\textbf{1.} $f(x) = e^{-ax}$, $a>0$. Then we come to the integral
\begin{equation}\label{example_int_1.1}
H_1(q) = \int\limits_0^{\infty} \!\!x e^{-ax} J_0(q x) \,dx  \;.
\end{equation}
On the other hand, $f^{(2m+1)}(0) = - a^{2m+1}$, and we find from
\eqref{H0_final} for $q > a$
\begin{equation}\label{example_sum_1.1}
H_1(q) = \frac{a}{q^3} \sum_{m=0}^\infty (-1)^m
\left(\frac{a}{2q}\right)^{\!2m} \frac{(2m+1)!}{(m!)^2} =
\frac{a}{(a^2+q^2)^{3/2}} \;,
\end{equation}
in accordance with equation \cite[2.12.8.4]{Prudnikov_vol_2}. To
calculate this series, we used formula
\cite[5.2.13.2]{Prudnikov_vol_1}.\\
\textbf{2.} $f(x) = x^{1+n} e^{-ax}$, $a>0$, integer $n \geqslant
0$. In such a case, $f^{(p)}(0) = 0$ for $0 \leqslant p \leqslant
n$, and $f^{(p)}(0) \neq 0$ for $p\geqslant n+1$.  We get
\cite[2.12.8.4]{Prudnikov_vol_2}
\begin{align}\label{example_int_2}
H_2(q) &= \int\limits_0^{\infty} \!\!x^{2+n} e^{-ax} J_0(q x) \,dx
\nonumber \\
&= \frac{1}{a^{n+3}} \Gamma(n+3) \,{}_2F_1 \!\left( \frac{n}{2} +
\frac{3}{2}, \frac{n}{2} + 2; 1; - \frac{q^2}{a^2} \right) ,
\end{align}
where ${}_2F_1 (a,b;c; z)$ is the hypergeometric function
\cite{Bateman_vol_1}. It satisfies equation
\begin{align}\label{hyper_fun_relation_2}
&{}_2F_1 \!\left( \frac{n}{2} + \frac{3}{2}, \frac{n}{2} + 2; 1; -
\frac{q^2}{a^2} \right) = \sqrt{\pi} \Bigg[ \frac{1}{\Gamma
\left((n+4)/2 \right)\Gamma \left(-(n+1)/2 \right)}
\nonumber \\
&\times \left( \frac{a}{q}\right)^{\!(n+3)} \!\!{}_2F_1 \!\left(
\frac{n}{2} + \frac{3}{2}, \frac{n}{2} + \frac{3}{2}; \frac{1}{2};
-\frac{a^2}{q^2} \right)- \frac{2}{\Gamma \left((n+3)/2
\right)\Gamma \left(-(n+2)/2 \right)}
\nonumber \\
&\times\left( \frac{a}{q}\right)^{\!(n+4)} \!\!{}_2F_1 \!\left(
\frac{n}{2} + 2, \frac{n}{2} + 2; \frac{3}{2}; -\frac{a^2}{q^2}
\right) \Bigg] .
\end{align}

First, take $n=2p$, $p=0, 1,2, \ldots$, then the r.h.s. of
\eqref{example_int_2} appears to be
\begin{align}\label{example_int_2_even_n}
H_2(q) &= \frac{(-1)^{p+1}}{2\pi}  \, [\Gamma(p+3/2)]^2 \left(
\frac{2}{q} \right)^{\!2p+3} {}_2F_1 \!\left( p + \frac{3}{2}, p +
\frac{3}{2}; \frac{1}{2}; -\frac{a^2}{q^2} \right) \;.
\end{align}
On the other hand, we have
\begin{equation}\label{derivatives_even_n}
f^{(2m+1)}(0) = \left\{
  \begin{array}{ll}
    0 , &  m<p \;, \\
    a^{2m-2p} \displaystyle{\frac{\Gamma(2m+2)}{\Gamma(2m-2p+1)}} \;, & m \geqslant
    p \;.
  \end{array}
\right.
\end{equation}
Then we obtain from \eqref{H0_final}
\begin{align}\label{example_sum_2_even_n_1}
H_2(q) &= \frac{(-1)^{p+1}}{2\pi} \, [\Gamma(p+3/2)]^2 \left(
\frac{2}{q} \right)^{\!2p+3}
\nonumber \\
&\times \sum_{k=0}^{\infty} \, \frac{(p+3/2)_k (p+3/2)_k}{(1/2)_k \,
k!} \left(- \frac{a^2}{q^2} \right)^{\!k} ,
\end{align}
where $(a)_n = \Gamma(a+n)/\Gamma(n)$ is the Pochhammer
symbol~\cite{Bateman_vol_1}. We see that
\eqref{example_sum_2_even_n_1} coincides with
\eqref{example_int_2_even_n}.

Analogously, we analyze the case $n=2p+1$, $p=0,1,2, \ldots,$ for
which
\begin{equation}\label{derivatives_odd_n}
f^{(2m+1)}(0) = \left\{
  \begin{array}{ll}
    0 , &  m<p+1 \;, \\
    -a^{2m-2p-1} \displaystyle{\frac{\Gamma(2m+2)}{\Gamma(2m-2p)}} \;, & m \geqslant
    p+1 \;.
  \end{array}
\right.
\end{equation}
It is easy to demonstrate that both eqs.~\eqref{example_int_2},
\eqref{hyper_fun_relation_2} and our formula \eqref{H0_final} give
the same expression
\begin{equation}\label{example_int_1.2_odd_n}
H_2(q) = \frac{(-1)^{p+1}}{2\pi} \, [\Gamma(p+5/2)]^2 a \left(
\frac{2}{q} \right)^{\!2p+5} {}_2F_1 \!\left( p + \frac{5}{2}, p +
\frac{5}{2}; \frac{3}{2}; -\frac{a^2}{q^2} \right) .
\end{equation}
\\
\textbf{3.} $f(x) = e^{-ax} I_0(c x)$, $a > c > 0$, then we must
consider the integral
\begin{equation}\label{example_int_3}
H_3(q) = \int\limits_0^{\infty} \!\!x e^{-ax} I_0(c x) J_0(q x) \,dx
\;,
\end{equation}
where $I_\nu(z)$ is the modified Bessel function of the first kind
\cite{Bateman_vol_2}. According to
\cite[2.15.16.2]{Prudnikov_vol_2}, this integral is equal to
\begin{equation}\label{example_I_3}
H_3(q) = \frac{2a}{\pi} \left(\frac{k}{q c}\right)^{\!\!3/2} \!\!(1
- k^2)^{3/4} \left[ 2\mathrm{E}(k) - \mathrm{K}(k) \right] .
\end{equation}
Here $\mathrm{K}(k)$ ($\mathrm{E}(k)$) is the complete elliptic
integral of the first (second) kind \cite{Bateman_vol_3}, and
\begin{equation}\label{k}
k = \frac{1}{\sqrt{2}} \left[
  1 - \frac{q^2 + a^2 - c^2}{\sqrt{(q^2 + a^2 - c^2)^2 + 4 q^2 c^2}}
  \right]^{1/2} .
\end{equation}
Note that $k \simeq c/q$ at large $q$. The expression in the r.h.s.
of \eqref{example_I_3} has the following asymptotics
\begin{align}\label{example_I_3_asym}
H_3(q)\Big|_{q \rightarrow \infty} &= \frac{a}{q^3} \Bigg[ 1
-\frac{3(2 a^2 + 3 c^2)}{4 q^2} + \frac{15 (8 a^4 + 40 a^2 c^2 + 15
c^4)}{64 q^4} \Bigg] + \mathrm{O}(q^{-9}) \;.
\end{align}

At the same time, we find that%
\footnote{For a special case $a=c$, this sum is given by formula
\cite[4.2.7.37]{Prudnikov_vol_1}.}
\begin{equation}\label{deriv_example_3}
f^{(n)}(0) = (-1)^{n} a^{n}\sum_{k=0}^{[n/2]} \binom{n}{2k}
\binom{2k}{k} \!\left( \frac{c^2}{4a^2} \right)^{\!\!k} ,
\end{equation}
where $\binom{n}{m}$ is the binomial coefficient, and $[z]$ means
the integer part of $z$. For odd-order derivatives, we obtain
\begin{equation}\label{odd_deriv_example_3}
f^{(2m+1)}(0) = - a^{2m+1}\sum_{k=0}^{m} \binom{2m+1}{2k}
\binom{2k}{k} \!\left( \frac{c^2}{4a^2} \right)^{\!\!k} .
\end{equation}
In particular,
\begin{align}\label{example_I_3_der}
f^{(1)}(0) &= - a \;,
\nonumber \\
f^{(3)}(0) &= - a\left[ a^2 + \frac{3 c^2}{2} \right] ,
\nonumber \\
f^{(5)}(0) &= - a\left[ a^4 + 5 a^2 c^2 + \frac{15 c^4}{8} \right] .
\end{align}
A compact analytical expression for $H_4(q)$ valid for $q > a+c$ is
calculated in Appendix~A to be
\begin{equation}\label{example_I_3_final}
H_3(q) = \frac{a}{q^3} \, F_4 \!\left( \frac{3}{2}, \frac{3}{2},
\frac{3}{2}, 1, -\frac{a^2}{q^2}, -\frac{c^2}{q^2} \right),
\end{equation}
where $F_4(\alpha, \beta, \gamma, \gamma', x, y)$ is the
hypergeometric series of two variables. It is the uniformly and
absolutely convergence series for $\sqrt{|x|} + \sqrt{|y|} < 1$
\cite{Bateman_vol_1}. Either from eqs.~\eqref{H0_final},
\eqref{odd_deriv_example_3} or from expansion
\eqref{example_I_3_final} we reproduce the same asymptotic expansion
\eqref{example_I_3_asym}.

As a byproduct, we obtained the new analytical expression
\eqref{example_I_3_final} for the improper integral
\eqref{example_int_3}, which is evidently much more appropriate for
evaluating its asymptotics at $q \gg a,c$ than the tabulated
expression \eqref{example_I_3}.

Let us underline that in all examples considered above, the series
\eqref{H0_final} \emph{converges uniformly and absolutely} for
$q>a$, and, consequently, it is the exact representation of integral
\eqref{H0_def}.

However, the series \eqref{H0_final} may diverge if $f^{(2m+1)}(0)$
rises rather quickly as $m \rightarrow\infty$ and violates
asymptotic condition \eqref{f_growth} of theorem~1. Suppose that
$f^{(n)}(0) = \mathrm{O}(a^n \Gamma(n+b))$, with $a,b>0$, as $n$
tends to infinity. In such a case, the series in
eq.~\eqref{H0_final} should be considered as \emph{an asymptotic
series} of $\mathcal{H}_0[f(x)](q)$, as $q \rightarrow \infty$ (see
also \cite{Willis}-\cite{MacKinnon}). Let us illustrate this point
by the following
example.\\
\textbf{4.} $f(x) = (x+a)^{-1}$, $a>0$. We have
\cite[2.12.3.6]{Prudnikov_vol_2}
\begin{equation}\label{example_int_4}
H_5(q) = \int\limits_0^{\infty} \!\!\frac{x}{x+a} \,J_0(q x) \,dx =
\frac{1}{q} - \frac{\pi a}{2} \,[\,\mathrm{\mathbf{H}}_0(a q) -
Y_0(a q)\,] \;,
\end{equation}
where $\mathrm{\mathbf{H}}_\nu(z)$ is the Struve function, and
$Y_\nu(z)$ is the Bessel function of the second kind
\cite{Bateman_vol_2}. We can use the asymptotic formula
\cite[10.42(2)]{Watson}
\begin{align}\label{Struve_Bessel_asymp} \left[
\,\mathrm{\mathbf{H}}_\nu(z) - Y_\nu (z) \, \right]
&\stackrel{\mathrm{as}}{=} \frac{1}{\pi} \sum_{m=0}^{p} \frac{\Gamma
(m + 1/2)}{\Gamma (\nu + 1/2 - m)} \left( \frac{z}{2}
\right)^{\!-(2m-\nu+1)}
\nonumber \\
&+ \mathrm{O}(z^{\!-(2p -\nu +3)}) \;, \quad |z| \rightarrow \infty,
\quad |\arg z| < \pi\;.
\end{align}
Then we find the asymptotics of integral \eqref{example_int_4}
\begin{equation}\label{example_int_4_asymp_1}
H_5(q) \stackrel{\mathrm{as}}{=} \frac{4}{\pi a^2 q^3} \left[
\sum_{m=0}^{p-1} (-1)^m\,2^{2m} \,[\Gamma (m + 3/2)]^2 (a q)^{\!-2m}
+ \mathrm{O}(q^{-2p}) \right] .
\end{equation}

On the other hand, $f^{(2m+1)}(0) = -a^{\!-2m-2} \Gamma(2m+2)$, and
we get from \eqref{H0_final}
\begin{align}\label{example_sum_4_2}
H_5(q) &\stackrel{\mathrm{as}}{=} \frac{4}{\pi a^2 q^3}
\sum_{m=0}^\infty (-1)^m\,2^{2m} \,[\Gamma (m + 3/2)]^2 (a
q)^{\!-2m}
\nonumber \\
&= \frac{1}{a^2 q^3} \,{}_3F_0 \!\left(\frac{3}{2}, \frac{3}{2},
1\,; -\frac{4}{(q a)^2}  \right) ,
\end{align}
which is just series \eqref{example_int_4_asymp_1}. Here
${}_{p}F_q(\alpha_1, \ldots \alpha_p ;\beta_1, \ldots \beta_q;z)$ is
the generalized hypergeometric series \cite{Bateman_vol_1}.

This series \eqref{example_sum_4_2} diverges, because $[\Gamma (m +
3/2)]^2 \simeq m^{2m} e^{-2m}$ as $m \rightarrow \infty$. Let us
define
\begin{align}\label{finite_sum}
R_n(q) &= \sum_{k=n+1}^{\infty} (-1)^m \,2^{2m} \,[\Gamma (m +
3/2)]^2
(a q)^{\!-2m} \nonumber \\
&= (-1)^n \,2^{2n+2} \,(a q)^{\!-2n-2} \sum_{p=0}^{\infty} (-1)^p
\,2^{2p} \,[\Gamma (n + p + 5/2)]^2 (a q)^{\!-2p} \;.
\end{align}
Since
\begin{align}\label{asymp_series}
\lim_{q \rightarrow \infty} q^{2n} R_n(q) &= 0 \;, \quad \mathrm{for
\ fixed \ } n \;, \nonumber \\
\lim_{n \rightarrow \infty} q^{2n} R_n(q) &= \infty \;, \ \ \!
\mathrm{for \ fixed \ } q \;,
\end{align}
we conclude that \eqref{example_sum_4_2} is by definition the
asymptotic power series of the integral $H_5(q)$ as $q \rightarrow
\infty$.

In all examples presented above, tabulated integrals were being
considered which were defined in terms of algebraic or special
functions. However, in some cases finding asymptotics of a tabulated
expression may need additional calculations, while formula
\eqref{H0_final} provides us with a desired asymptotic expansion
immediately.

We can illustrate this statement by the following example.\\
\textbf{5.} $f(x) = (x+a)^{-2}$, $a>0$. It means that we consider
the integral
\begin{equation}\label{example_int_5}
H_5(q) = \int\limits_0^{\infty} \!\!\frac{x}{(x+a)^2} \,J_0(q x)
\,dx \;.
\end{equation}
It is a tabulated integral \cite[2.12.3.9]{Prudnikov_vol_2}. But
finding its asymptotics with the use of tabulated expression faces
technical difficulties, since one should start from the integral
\begin{equation}\label{example_int_5_eps}
H_5(q, \varepsilon) = \int\limits_0^{\infty}
\!\!\frac{x}{(x+a)^{2+\varepsilon}} \,J_0(q x) \,dx \;,
\end{equation}
calculate its asymptotics for $\varepsilon \neq 0$, and only then
take the limit $\varepsilon\rightarrow0$. The reason is that the
tabulated expression for integral \eqref{example_int_5_eps} has
three terms two of which having simple poles at $\varepsilon = 0$
\cite{Prudnikov_vol_2}.
That is why, one cannot put $\varepsilon=0$ from the very beginning.%
\footnote{In a sum of three terms, the poles cancel each other.}

On the other hand, formula \eqref{H0_final} gives us the asymptotic
expansion of $H_5(q)$ \eqref{example_int_5} as $q \rightarrow
\infty$, without additional calculations. Indeed, from relation
$f^{(2m+1)}(0) = -a^{\!-2m-3} \Gamma(2m+3)$ we immediately find the
desired expansion
\begin{align}\label{example_sum_1.5}
H_5(q) &\stackrel{\mathrm{as}}{=} \frac{8}{\pi (q a)^3}
\sum_{m=0}^\infty (-1)^m\,2^{2m} \,[\Gamma (m + 3/2)]^2 (m+1) (a
q)^{\!-2m}
\nonumber \\
&= \frac{2}{(q a)^3} \,{}_3F_0 \!\left(\frac{3}{2}, \frac{3}{2},
2\,; -\frac{4}{q^2 a^2}  \right) .
\end{align}
Note that this series can be obtained by differentiating series
\eqref{example_sum_4_2} with respect to the parameter $a$, since
$H_5 = -(\partial/\partial a)H_4$.

\textsc{Remark~1}. Theorem~1 remains true if $f^{(2m)}(0) = 0$,
$m=0,1,2, \ldots$. This statement can be illustrated, for example,
by the Hankel transform of the function $f(x) = J_1(a x) J_0(b x)$,
$a,b > 0$.

\textsc{Remark~2}. The condition \textbf{(2)} of theorem~1 might be
weakened, if one imposes
it on the function $\bar{f}(x) = x f(x)$, $\bar{f}(-x) \neq -\bar{f}(x)$.%
\footnote{It corresponds to an alternative definition of the Hankel
transform: $\bar{\mathcal{H}}_\nu [f(x)](q) = \int_0^{\infty} \!\!
f(x) J_\nu(q x) \,dx$.}
Correspondingly, the condition \textbf{(3)} must be: $x^{-1/2}
\bar{f}(x) \rightarrow 0$, as $x \rightarrow \infty$. In such a
case, the Hankel transform has the form
\begin{equation}\label{H0_m}
\mathcal{H}_0 [f(x)](q) = \frac{1}{q} \sum_{m=0}^\infty (-1)^{m}
\,\frac{\Gamma (2m+1)}{\Gamma^2 (m+1)} \,\bar{f}^{(2m)}(0) \,
(2q)^{-2m} \;.
\end{equation}
As an example, let us consider the function $f(x) = e^{-ax}\!/x$,
then $\bar{f}^{(2m)}(0) = a^{2m}$. Using \eqref{H0_m} and equation
\cite[5.2.13.1]{Prudnikov_vol_1}, we find
\begin{equation}\label{modified_condition}
\mathcal{H}_0 [e^{-ax}\!/x](q) = \frac{1}{\sqrt{a^2 + q^2}} \;,
\end{equation}
in agreement with \cite[2.12.8.3]{Bateman_vol_2}.

\section{Hankel transform of order one and two} %
\label{sec:one}

\textsc{Theorem} 2. The first-order Hankel transform
\begin{equation}\label{H1_def}
\mathcal{H}_1 [f(x)](q) = \int\limits_0^{\infty} \!\!x f(x) J_1(q
x)\,dx
\end{equation}
can be presented by absolutely and uniformly convergent series for
$q>a>0$
\begin{align}\label{H1_final}
\mathcal{H}_1 [f(x)](q) &= \frac{1}{q^2} \sum_{m=0}^\infty (-1)^{m}
\frac{\Gamma (2m+2)}{\Gamma^2 (m+1)} f^{(2m)}(0) \, (2q)^{-2m} \;,
\nonumber \\
&= \frac{2}{q^2} \sum_{k=0}^\infty \,\frac{\Gamma (k/2 +
3/2)}{\Gamma (1/2 - k/2) \,\Gamma(k+1)} \,f^{(k)}(0) \left(
\frac{q}{2} \right)^{\!\!-k}  \;.
\end{align}
provided that $f(x)$ and its derivatives  at $x=0$ satisfy
conditions \textbf{(2)}-\textbf{(4)} of theorem~1, and the condition
\textbf{(1)} is replaced by the condition $f(-x) \neq - f(x)$.

\textsc{Proof}. The proof of this theorem can be done by analogy
with the proof of theorem~1. Another simple way is to use the
following relation \cite{Bateman_vol_2}
\begin{equation}\label{Bessel0_relation}
J_0(q x) = \frac{1}{q} \left[ \frac{1}{x} \,J_1(q x) +  \frac{d}{dx}
J_1(q x) \right] .
\end{equation}
Then we have%
\footnote{We used the fact that $xf(x)J_1(q x)|_{x=0} = xf(x)J_1(q
x)|_{x=\infty} =0$.}
\begin{equation}\label{J0_vsJ1}
\int\limits_0^{\infty} \!x f^{(1)}(x) J_1(q x) \,dx = - q \!\!
\int\limits_0^{\infty} \!x f(x) J_0(q x) \,dx \;.
\end{equation}
Formula~\eqref{H1_final} follows immediately from
eqs.~\eqref{J0_vsJ1} and \eqref{H0_final}. \textsc{Q.E.D.}

\textsc{Remark~3}. The Hankel transform $\mathcal{H}_{-1} [f(x)](q)$
is given by formula \eqref{H1_final} taken with the opposite sign,
since $J_{-1}(x) = -J_1(x)$.

Let us consider two examples.\\
\textbf{6.}  $f(x) = e^{-a x}$, $a>0$, then
\begin{equation}\label{example_int_6}
H_6(q) = \int\limits_0^{\infty} \!\!x e^{-a x} J_1(q x) \,dx  \;.
\end{equation}
Taking into account that $f^{(2m)}(0) = a^{2m}$, we immediately
obtain from \eqref{H1_final} that
\begin{equation}\label{example_final_6}
H_6(q) = \frac{q}{(a^2+q^2)^{3/2}}
\end{equation}
for $q>a$, in agreement with \cite[2.12.8.4]{Prudnikov_vol_2}.\\
\textbf{7.} $f(x) = e^{-ax} I_0(c x)$, $a > c > 0$, then we deal
with the integral
\begin{equation}\label{example_int_7}
H_7(q) = \int\limits_0^{\infty} \!\!x e^{-ax} I_0(c x) J_1(q x) \,dx
\;,
\end{equation}
The even-order derivatives of $f(x)$ are given by the formula (see
eq. \eqref{deriv_example_3})
\begin{equation}\label{odd_deriv_example_7}
f^{(2m)}(0) = a^{2m}\sum_{k=0}^{m} \binom{2m}{2k} \binom{2k}{k}
\!\left( \frac{c^2}{4a^2} \right)^{\!\!k} .
\end{equation}
Doing in the same manner as in considering example 3, we find
\begin{equation}\label{example_I_7_final}
H_7(q) = \frac{1}{q^2} \, F_4 \!\left( \frac{3}{2}, \frac{1}{2}, 1,
\frac{1}{2}, -\frac{a^2}{q^2}, -\frac{c^2}{q^2} \right) .
\end{equation}
This expression is much more appropriate for evaluating the
asymptotics of \eqref{example_int_7} at large $q$ than the
corresponding tabulated expression
\cite[2.15.16.2]{Prudnikov_vol_2}.

\textsc{Theorem} 3. The second-order Hankel transform
\begin{equation}\label{H2_def}
\mathcal{H}_2 [f(x)](q) = \int\limits_0^{\infty} \!\!x f(x) J_2(q
x)\,dx
\end{equation}
can be presented by absolutely and uniformly convergent series for
$q>a>0$
\begin{align}\label{H2_final}
\mathcal{H}_2 [f(x)](q) &= \frac{2}{q^2} f(0) + \frac{4}{\pi q^3}
\sum_{m=0}^\infty (-1)^{m} \,\frac{\Gamma(m+5/2)
\,\Gamma(m+1/2)}{\Gamma(2m+2)} \nonumber \\
&\times f^{(2m+1)}(0) \, \left( \frac{q}{2} \right)^{-2m}  \nonumber \\
&= \frac{2}{q^2} \sum_{k=0}^\infty \,\frac{\Gamma (k/2 + 2)}{\Gamma
(1 - k/2) \,\Gamma(k+1)} \,f^{(k)}(0) \left( \frac{q}{2}
\right)^{\!\!-k} ,
\end{align}
provided that $f(x)$ and its derivatives  at $x=0$ satisfy all the
conditions of theorem~1.

\textsc{Proof}. Observe that
\begin{equation}\label{residue}
\frac{2}{q^2} \,f(0) = \pi\mathrm{Res}[F(-s) \left( \frac{q}{2}
\right)^{2s} , 1/2] \;,
\end{equation}
where
\begin{equation}\label{F}
F(s) = \frac{4}{\pi q^3} \,\frac{\Gamma(s+5/2)
\,\Gamma(s+1/2)}{\Gamma(2s+2)} \,f^{(2s+1)}(0) \;.
\end{equation}
Note that $f(0) \equiv f^{(0)}(0)$. The function $F(z)$ satisfies
all the condition of the Ramanujan master theorem. According to the
corollary of this theorem (see section \ref{sec:Ramanujan}), we have
\begin{align}\label{S2}
S_2(q) &= \frac{2}{q^2} f(0) + \frac{4}{\pi q^3} \sum_{m=0}^\infty
(-1)^{m} \,\frac{\Gamma(m+5/2)
\,\Gamma(m+1/2)}{\Gamma(2m+2)} \,f^{(2m+1)}(0) \!\left( \frac{q}{2} \right)^{-2m}  \nonumber \\
&= \frac{i}{4\pi} \int\limits_{c-i\infty}^{c+i\infty}
\frac{\Gamma(s+5/2) \,\Gamma(s+1/2)}{\sin(\pi s)\Gamma(2s+1)}
\,g^{(2s+1)}(0) \!\left( \frac{q}{2} \right)^{-2s-3}  ds \;,
\end{align}
where $-\delta < c < -3/4$. The remaining part of the proof proceeds
as the proof of theorem~1. After change of variables $2s+1 = s'$ and
simple calculations, we obtain
\begin{equation}\label{S2_contour_median}
S_2(q) = \frac{1}{4i\pi} \int\limits_{c'-i\infty}^{c'+i\infty}
\frac{ \Gamma(-s') \Gamma(s'/2 + 2)}{\Gamma(-s'/2 + 1)}
\,g^{(s')}(0) \,\left( \frac{q^2}{4} \right)^{\!\!-(s'/2+1)} ds' \;,
\end{equation}
where $-1 < c' < -1/2$. Using formula \cite[2.12.2.2]{Bateman_vol_2}
\begin{equation}\label{Bessel_nu_integral}
\frac{\Gamma(s/2 + 2)}{\Gamma(-s/2 + 1)} \, \left( \frac{q^2}{4}
\right)^{\!\!-(s/2+1)} = 2 \!\int\limits_0^\infty \!\!x^{s+1} J_2 (x
q) \,dx \;,
\end{equation}
valid for $-4 < \mathfrak{R}(s) < -1/2$, and the Cauchy residue
theorem, we come to \eqref{H2_def}. Q.E.D.

The conditions of theorems~1-2 guarantee that the r.h.s of
eqs.~\eqref{H1_final}, \eqref{H2_final} are \emph{the uniformly
convergent series}, not only asymptotic series.

\section{Hankel transform of an arbitrary integer order} %
\label{sec:nu}

\textsc{Theorem 4}. The Hankel transform of order $n \geqslant -1$
can be presented by absolutely and uniformly convergent series for
$q>a>0$
\begin{equation}\label{Hnu_final}
\mathcal{H}_n [f(x)](q) = \frac{2}{q^2} \sum_{k=0}^\infty
\,\frac{\Gamma (k/2 + n/2 + 1)}{\Gamma (n/2 - k/2) \,\Gamma(k+1)}
\,f^{(k)}(0) \left( \frac{q}{2} \right)^{\!\!-k} ,
\end{equation}
provided that $f(x)$ and its derivatives at $x=0$ satisfy all the
conditions of theorem~1, except for the condition \textbf{(1)} which
is replaced by the following condition
\begin{equation}\label{condition_1_mod}
f(-x) \neq  \left\{
                   \begin{array}{ll}
                   + f(x)  \;, &  n \mathrm{\ even} \;, \\
                   - f(x) \;, & n \mathrm{\ odd} \;.
                   \end{array}
             \right.
\end{equation}

\textsc{Proof}. We apply the strong induction which is a
strengthening of the basic mathematical induction \cite{math_ind}.
Let $P(n)$ be the statement of theorem~4 for the Hankel transform of
order $n$. According to theorems~2 and 3, P(1) and $P(2)$ hold.
Suppose that for every $m+1 > 2$, $P(k)$ is true for all $k < m+1$.
Let us show that $P(m+1)$ is also true.

Form known relations between the Bessel functions one can easily
find that for $\nu > -1$%
\footnote{Taking into account that $xf(x)J_\nu(q x)|_{x=0} =
xf(x)J_\nu(q x)|_{x=\infty} =0$.}
\begin{align}\label{Bessel_relation}
(\nu-1) \int\limits_0^{\infty} \!x f(x) J_{\nu+1}(q x)\,dx &=
(\nu+1) \int\limits_0^{\infty} \!x f(x) J_{\nu-1}(q
x)\,dx \nonumber \\
&+ \frac{2\nu}{q} \int\limits_0^{\infty} \!x f^{(1)}(x) J_{\nu}(q
x)\,dx \;.
\end{align}
Since $P(m)$ and $P(m-1)$ are assumed to be true, we have
\begin{align}\label{Hm-1}
(m+1)\mathcal{H}_{m-1} [f(x)](q) &= \frac{2(m+1)}{q^2}
\sum_{k=0}^\infty \,\frac{\Gamma (k/2 + m/2 + 1/2)}{\Gamma (m/2 -
k/2 - 1/2) \,\Gamma(k+1)} \nonumber \\
&\times f^{(k)}(0) \left( \frac{q}{2} \right)^{\!\!-k} ,
\end{align}
and
\begin{align}\label{Hm}
\frac{2m}{q} \mathcal{H}_m [f^{(1)}(x)](q) &= \frac{2m}{q^2}
\sum_{k=-1}^\infty \,\frac{\Gamma (k/2 + m/2 + 1)}{\Gamma (m/2 -
k/2) \,\Gamma(k+1)} \,f^{(k+1)}(0) \left( \frac{q}{2}
\right)^{\!\!-(k+1)} \nonumber \\
&= \frac{2m}{q^2} \sum_{k'=0}^\infty \,\frac{\Gamma (k'/2 + m/2 +
1/2)}{\Gamma (m/2 - k'/2 + 1/2) \,\Gamma(k')} \,f^{(k')}(0) \left(
\frac{q}{2} \right)^{\!\!-(k')} .
\end{align}
Then we find from \eqref{Bessel_relation}-\eqref{Hm}, using relation
$(m+1)(m/2 - k/2 - 1/2) + m k = (m-1)(m/2 + k/2 + 1/2)$,
\begin{align}\label{Pm+1}
\mathcal{H}_{m+1} [f(x)](q) &= \frac{1}{m-1} \left[
(m+1)\mathcal{H}_{m-1} [f(x)](q) +
\frac{2m}{q} \mathcal{H}_m [f^{(1)}(x)](q) \right] \nonumber \\
&= \frac{2}{q^2} \sum_{k=0}^\infty \,\frac{\Gamma (k/2 + m/2 +
3/2)}{\Gamma (m/2 - k/2 + 1/2 ) \,\Gamma(k+1)} \,f^{(k)}(0) \left(
\frac{q}{2} \right)^{\!\!-k} \nonumber \\
&= P(m+1) \;.
\end{align}
It means that the inductive step holds. Thus, the statement of
theorem~4 is true for every integer $n \geqslant -1$. Q.E.D.

A formula analogous to \eqref{Hnu_final} was obtained in
\cite[(5.1)]{MacKinnon}, but as \emph{asymptotic expansion} of
$\mathcal{H}_\nu [f(x)](q)$ as $q$ tends to infinity. In our case,
the r.h.s. of \eqref{Hnu_final} is the uniformly and absolutely
convergent power series  for every $q>a>0$. Thus, it is the exact
expansion of the Hankel transform $\mathcal{H}_n [f(x)](q)$.

\textbf{8.} As an example, we consider the integral%
\footnote{Two particular cases with $n=0$ and $n=1$ were considered
above (see eqs.~\eqref{example_int_1.1}, \eqref{example_int_7}).}
\begin{equation}\label{example_int_8}
H_8(q) = \int\limits_0^{\infty} \!\!x e^{-ax} J_n(q x) \,dx  \;,
\quad n\geqslant 2 \;.
\end{equation}
We find from \eqref{Hnu_final}, using $f^{(n)}(0) = (-1)^n a^{n}$,
\begin{align}\label{H8_1}
H_8(q) &= \frac{2}{q^2} \,\sum_{k=0}^\infty \frac{1}{k!} \,
\frac{\Gamma (k/2 + n/2
 + 1))}{\Gamma (n/2 - k/2)}  \left( -\frac{2a}{q} \right)^{\!\!k}
\nonumber \\
&= \frac{1}{q^2} \left( n + \frac{\partial}{\partial \ln x} \right)
\sum_{k=0}^\infty \frac{1}{k!} \, \frac{\Gamma (k/2 +
n/2)}{\Gamma(n/2 - k/2)} \,(-2\,x)^k \;,
\end{align}
where $x = a/q$. The series on the second line of eq.~\eqref{H8_1}
is equal to \cite[5.2.14.29]{Prudnikov_vol_1}
\begin{equation}\label{S}
\sum_{k=0}^\infty \frac{1}{k!} \, \frac{\Gamma (k/2 + n/2)}{\Gamma
(n/2 - k/2)}  \,(-2\,x)^k = \frac{(\sqrt{1+x^2} - x)^{n-1}}{\sqrt{1
+ x^2}} \;.
\end{equation}
This series converges for $|x| < 1$. As a result, we come to the
expression
\begin{align}\label{H8_final}
H_8(q) &= \frac{1}{q^2} \left( n + \frac{\partial}{\partial \ln x}
\right) \!\frac{(\sqrt{1 + x^2} - x)^{n-1}}{\sqrt{1 + x^2}} \nonumber \\
&= -\frac{1}{q^2} \, \frac{\partial}{\partial x} \frac{(\sqrt{1 +
x^2} - x)^{n}}{\sqrt{1 + x^2}} \;.
\end{align}
Formula \cite[2.12.8.4]{Prudnikov_vol_2} (with $\alpha = 2$) gives
the same result for $H_8(q)$ \eqref{H8_final}. It demonstrates us
that series \eqref{Hnu_final} is really the exact representation of
the Hankel transform \eqref{Hnu_def} for all $q>a>0$.

As one can see from \eqref{Hnu_final}, in the large-$q$ limit
$\mathcal{H}_n[f(x)](q)$ is $\mathrm{O}(q^{-3})$ for $n=0$ and
$\mathrm{O}(q^{-2})$ for $n = -1, 1, 2, \ldots$.


\section{Discussions and conclusions} %

We have studied the Hankel transform \eqref{Hnu_def} for $\nu = n
\geqslant -1$ and positive $q$, imposing several conditions on
$f(x)$ and its derivatives at $x=0$. In particular, we demand that
its $s$-th derivative (up to a sign) $g^{(s)}(0)$ is a regular
function defined on a half-plane $\mathfrak{R}(s) > - \delta$ for
some $3/4 < \delta < 1$, and that it grows as $a^{\mathfrak{R}(s)}$
when $s \rightarrow \infty$ (see the conditions of theorem~1).

It is shown that the Hankel transform can be presented by the
absolutely and uniformly convergent series in reciprocal powers of
$q$ for $q > a > 0$ (see eqs.~\eqref{H0_final}, \eqref{H1_final},
\eqref{H2_final}, and \eqref{Hnu_final}). If one of the conditions
imposed on $f(s)$ is violated, series \eqref{Hnu_final} becomes
asymptotic series as $q \rightarrow \infty$.

The validity of our theorems has been illustrated by eight examples.
As a byproduct, we have obtained the new analytical expression
\eqref{example_I_3_final} for the definite integral of the Bessel
functions $J_0(z)$ and $I_0(z)$ \eqref{example_int_3}, which is much
more appropriate for evaluating asymptotics of \eqref{example_int_3}
at large $q$ than the tabulated expression \eqref{example_I_3}. The
same must be said about our expression \eqref{example_I_7_final}. In
some cases finding asymptotics of a tabulated integral may need
complicated calculations, while our exact expansion provides us with
its asymptotics immediately, as it is shown for integral
\eqref{example_int_5}.

Previously, in a number of papers similar formulas were derived
\cite{Willis}-\cite{Soni}. However, it is necessary to underline
that they are in fact \emph{asymptotic expansions} valid only when
$q$ tends to infinity.

In conclusion, one result of \cite{Hsu} is worth discussing. The
author of \cite{Hsu} studied the integral
\begin{equation}\label{Hsu_int_def}
I(\lambda) = \int\limits_0^\infty \!\Phi(\lambda t) f(t) \,dt \;,
\end{equation}
where $\lambda$ is a positive parameter, and $\Phi(t)$ is assumed to
have a Laplace transform $\Psi(s)$. $\Phi(t)$ can be one of familiar
functions, e.g. the Bessel function $J_\nu(t)$. Let $f(t) =
\sum_{n=0}^\infty c_n t^n$ be an entire (i.e. analytic at all finite
points of the complex plane $\mathbb{C}$) function such that both
$\Phi(\lambda t) f(t)$ and
\begin{equation}\label{integrand_abs}
|\Phi(\lambda t)| \sum_{n=0}^\infty |c_n| t^n
\end{equation}
are integrable over $[0,\infty)$. Then theorem~1 in \cite{Hsu} says
that
\begin{equation}\label{Hsu_int_final}
I(\lambda) =   \sum_{n=0}^\infty \frac{(-1)^n}{n!} \, \Psi^{(n)}(0)
\,f^{(n)}(0) \,\lambda^{-n-1} \;.
\end{equation}

However, \eqref{integrand_abs} \emph{isn't integrable} for most
integrals investigated in the present paper.%
\footnote{See integrals \eqref{example_int_1.1},
\eqref{example_int_2}, \eqref{example_int_3}, \eqref{example_int_6}
and \eqref{example_int_8}.}
Consider, for instance, integral $H_1$ \eqref{example_int_1.1}. In
notations of \cite{Hsu}, $\Phi(t) = J_0(t)$, and
\begin{equation}\label{Hsu_H1}
f(t) = te^{-at} = \sum_{n=0}^\infty c_n t^n \;,
\end{equation}
where $c_n = (-a)^{n-1}/\Gamma(n)$. We get
\begin{equation}\label{H1_int_def}
|\Phi(\lambda t)| \sum_{n=0}^\infty |c_n| t^n = |J_0(\lambda t)|
\,|t \,e^{at}| \;.
\end{equation}
It is evident that function \eqref{H1_int_def} isn't integrable over
$[0,\infty)$ for $\lambda, a>0$.%
\footnote{The same is true for $\Phi(t) = J_0(t)$, $f(t) = t\cos t,
t\sin t, tJ_\nu(t)$, \emph{etc}.}
So, one of conditions of theorem~1 in \cite{Hsu} is violated. On the
contrary, our theorem~1 works, and it results in the exact
expression for $H_1(q)$ \eqref{example_sum_1.1} .

All said above allows us to conclude that our formulas
\eqref{H0_final}, \eqref{H1_final}, \eqref{H2_final}, and
\eqref{Hnu_final} are the new results which can be used to obtain
analytic expressions for the Hankel transform
$\mathcal{H}_n[f(x)](q)$ with the positive parameter $q$.



\setcounter{equation}{0}
\renewcommand{\theequation}{A.\arabic{equation}}

\section*{Appendix A}
\label{app:A}

Here we calculate a series which defines integral
\eqref{example_int_3}. After using the relation
\begin{equation}\label{Gamma_relation}
\Gamma(2m+2) = 2^{2m+1} \Gamma(m+1) \Gamma(m+3/2)/ \sqrt{\pi} \;,
\end{equation}
we find from \eqref{H0_final} and \eqref{odd_deriv_example_3}
\begin{align}\label{I3_def}
&H_3(q) = \frac{2a}{\sqrt{\pi} q^3} \sum_{m=0}^\infty
\frac{\Gamma(m+3/2)}{\Gamma(m+1)}
\left(-\frac{a^2}{q^2}\right)^{\!\!m} \sum_{k=0}^m \binom{2m+1}{2k}
\binom{2k}{k} \!\left(\frac{c^2}{4a^2}\right)^{\!\!k}
\nonumber \\
&= \frac{2a}{\sqrt{\pi} q^3} \sum_{k=0}^\infty \frac{1}{(k!)^2}
\left(\frac{c^2}{4a^2}\right)^{\!\!k} \sum_{m=k}^\infty
\frac{\Gamma(m+3/2)\,\Gamma(2m+2)}{\Gamma(m+1)\,\Gamma(2m+2-2k)}
\left(-\frac{a^2}{q^2}\right)^{\!\!m} .
\end{align}
After putting $m=n+k$, we find
\begin{align}\label{I3_1}
H_3(q) &=\frac{2a}{\sqrt{\pi} q^3} \sum_{k=0}^\infty
\frac{1}{(k!)^2} \left(-\frac{c^2}{4q^2}\right)^{\!\!k}
\nonumber \\
&\times \sum_{n=0}^\infty
\frac{\Gamma(n+k+3/2)\,\Gamma(2n+2k+2)}{\Gamma(n+k+1)\,\Gamma(2n+2)}
\left(-\frac{a^2}{q^2}\right)^{\!\!n} .
\end{align}
Using relation
\begin{equation}\label{gamma_fun_1}
\frac{\Gamma(n+k+3/2)\,\Gamma(2n+2k+2)}{\Gamma(n+k+1)\,\Gamma(2n+2)}
= 2^{2k} \, \frac{[\Gamma(n+k+3/2)]^2}{\Gamma(n+3/2) \, n!} \;,
\end{equation}
we obtain ($q > a+c$)
\begin{align}\label{I3_2}
H_3(q) &= \frac{a}{q^3} \sum_{k=0}^\infty \frac{1}{(1)_k k!}
\left(-\frac{c^2}{q^2}\right)^{\!\!k} \sum_{n=0}^\infty
\frac{[(3/2)_{n+k}]^2}{(3/2)_n n!}
\left(-\frac{a^2}{q^2}\right)^{\!\!n}
\nonumber \\
&= \frac{a}{q^3} \, F_4 \!\left( \frac{3}{2}, \frac{3}{2},
\frac{3}{2}, 1, -\frac{a^2}{q^2}, -\frac{c^2}{q^2} \right) ,
\end{align}
where
\begin{equation}\label{hyper_geom_xy}
F_4(\alpha, \beta, \gamma, \gamma', x, y) = \sum_{n.m=0}^\infty
\frac{(\alpha)_{m+n} (\beta)_{m+n}}{(\gamma)_{m} (\gamma')_{n} \,m!
\,n!} \, x^m y^n
\end{equation}
is the hypergeometric series of two variables \cite{Bateman_vol_1}.
It is the uniformly and absolutely convergence series for
$\sqrt{|x|} + \sqrt{|y|} < 1$.




\end{document}